\documentclass[11pt]{amsart}
\usepackage[active]{srcltx}
\usepackage{color}
\usepackage{graphicx}
\usepackage{caption}
\usepackage{subcaption}

\theoremstyle{plain}
\newtheorem{theorem}{Theorem}
\newtheorem{lemma}[theorem]{Lemma}
\newtheorem{proposition}[theorem]{Proposition}
\newtheorem{corollary}[theorem]{Corollary}

\theoremstyle{definition}
\newtheorem{definition}[theorem]{Definition} \theoremstyle{remark}

\begin{document}
 
\newcommand{\N}{{\mathbb N}}

\def\A{{\mathcal A}}
\def\S{{\mathcal S}}
\def\H{{\mathcal H}}
\def\M{{\mathcal M}}
\def\MM{{\mathcal M}}
\def\T{{\mathcal T}}
\def\I{{\mathcal I}}
\def\F{{\mathcal F}}
\def\J{{\mathcal J}}
\def\E{{\mathcal E}}
\def\P{{\mathcal P}}
\def\HH{{\mathcal H}}
\def\V{{\mathcal V}}
\def\B{{\mathcal B}}
\def\grad{\nabla}
\font\teneufm=eufm10
\font\seveneufm=eufm7
\font\fiveeufm=eufm5
\newfam\eufmfam
\textfont\eufmfam=\teneufm
\scriptfont\eufmfam=\seveneufm
\scriptscriptfont\eufmfam=\fiveeufm
\def\eufm#1{{\fam\eufmfam\relax#1}}
\def\Def{{\eufm D}}
\def\Man{{\eufm M}}

\newcommand{\average}{{\mathchoice {\kern1ex\vcenter{\hrule
height.4pt width 6pt depth0pt}
\kern-11pt} {\kern1ex\vcenter{\hrule height.4pt width 4.3pt
depth0pt} \kern-7pt} {} {} }}
\newcommand{\ave}{\average\int}

\title[Recent Results On Nonlinear Elliptic Free Boundary Problems]{Recent Results On Nonlinear Elliptic Free Boundary Problems}
\author{Fausto Ferrari}
\address{Dipartimento di Matematica dell'Universit\`a di Bologna, Piazza di Porta S. Donato, 5, 40126 Bologna, Italy.}
\email{\tt fausto.ferrari@unibo.it}
\author[Claudia Lederman]{Claudia Lederman}
\address{IMAS - CONICET and Departamento  de
Ma\-te\-m\'a\-ti\-ca, Facultad de Ciencias Exactas y Naturales,
Universidad de Buenos Aires, (1428) Buenos Aires, Argentina.}
\email{\tt  clederma@dm.uba.ar}
\author{Sandro Salsa}
\address{Dipartimento di Matematica del Politecnico di Milano, Piazza Leonardo da Vinci, 32, 20133 Milano, Italy.}
\email{\tt sandro.salsa@polimi.it }

\thanks{F.F. was partially supported by INDAM-GNAMPA project 2020, Metodi di viscosit\`a e applicazioni a pro\-ble\-mi non lineari con debole ellitticit\`a. C. L. was partially supported by GHAIA Horizon 2020 MCSA RISE 2017
programme grant 777822, CONICET PIP 11220150100032CO 2016-2019, UBACYT
20020150100154BA and ANPCyT PICT 2016-1022.}

\maketitle

\emph{To Alfio Quarteroni: a great mathematician, a great person, a great friend}

\begin{abstract}
In this paper we give an overview   of some recent and older results concerning free boundary problems governed by elliptic operators.
\end{abstract}

\section{Introduction}

In this survey paper we consider free boundary problems governed by elliptic
equations in which the state variable can assume two phases and the
condition across the free boundary is expressed by an energy balance
involving the fluxes from both sides (so called Bernoulli type problems).
Typical examples come from constraint minimization problems, flame
propagation as limits of singular perturbation problems with forcing term
(flame propagation), \cite{LW, T}, the Stokes and Prandtl-Batchelor
models in classical hydrodynamics, \cite{AFT, B1}, free transmission
problems, \cite{AT}.

Among the several concepts of solutions, we use the notion of viscosity
solution introduced by Luis Caffarelli in the seminal papers \cite{C1,C2, C3}%
, which seems to be the most appropriate to study optimal regularity of both
solutions and free boundaries in great generality. In those papers, in the
homogeneous case, Caffarelli developed a general strategy to attack the
existence and the optimal regularity of both the solution and the free
boundary based on a powerful monotonicity formula, already proved in \cite{ACF}, Harnack principles and
families of continuous perturbations.

A different approach to regularity, relying on Harnack principles and
linearization, has been introduced by De Silva in \cite{D}, for a one phase
model motivated by the classical Stokes problem in hydrodynamics, and
subsequently refined in \cite{DFS}, \cite{DFS1}, \cite{DFS2}, to cover a
broad spectrum of applications, in particular, for problems with distributed
sources.

The main common underlying idea of both techniques to obtain the regularity
of the free boundary, frequently written simply f.b. in the sequel, is to set up an iterative improvement of flatness
argument in a neighborhood of a point where one of the two phases enjoys a
non-degeneracy condition (e.g. a linear growth).

After the quoted seminal papers by Caffarelli, the theory can be developed
according to a well established paradigm:

\begin{itemize}
\item[a).] Existence and optimal regularity of solutions, e.g. viscosity or
variational solutions, or solutions obtained as a limit of singular
perturbations.

\item[b).] Weak regularity properties of the f.b., such as finite perimeter
and density properties for the positivity set.

\item[c).] Strong regularity properties of the f.b. For instance Lipschitz
or ``flat" free boundaries are $C^{1}$ or better.

\item[d).] Higher regularity: Schauder type estimates and analyticity for
both solution and f.b.
\end{itemize}

Also of great importance, we believe, is to have information on the
Hausdorff measure or dimension of the \emph{singular} (\emph{nonflat})
points of the free boundary. The best existing results in this
direction are obtained by De Philippis, Spolaor, Velichkov for the classical
Bernoulli problem in \cite{DSV}. Nothing is known in the nonhomogeneous case.%

\bigskip

In Sections \ref{sec2} to \ref{sec5} of this paper we focus on f.b. problems governed
by uniformly elliptic equation, while in Section \ref{sec6} we describe some
recent results on nonlinear operators with non-standard growth, see \cite{FL}.

\section{Free boundary problems governed by fully nonlinear operators}

\label{sec2}

In the last years the theory for two phase problems governed by uniformly
elliptic operators has reached a considerable level of completeness. Here we
mostly focus on problems governed by fully nonlinear equations with
distributed sources. Our class of free boundary problems and their viscosity
solutions can be stated as follows.

Let ${\mathrm{Sym}}_{n}$ denote the space of $n\times n$ symmetric matrices
and let $F:{\mathrm{Sym}}_{n}\rightarrow \mathbb{R}$ with $F\left( O\right)
=0$ and such that there exist constants $0<\lambda \leq \Lambda $ with 
\begin{equation*}
\lambda \Vert N\Vert \leq F(M+N)-F(M)\leq \Lambda \Vert N\Vert \qquad \text{
for every }M,N\in {\mathrm{Sym}}_{n}\text{ with }N\geq 0,
\end{equation*}%
where $\Vert M\Vert =\max_{|x|=1}|Mx|$ denotes the $(L^{2},L^{2})$-norm of
the matrix $M$. \ Observe that if $\mathcal{F}$ is an operator in our class
then for every $r>0$ 
\begin{equation}\label{operatorname}
F_{r}(M)=\frac{1}{r}F(rM)
\end{equation}%
is still an operator in the same class.\bigskip

Let $\Omega \subset {\mathbb{R}}^{n}$ be a bounded Lipschitz domain and $%
f_{1},f_{2}\in C(\Omega )\cap L^{\infty }(\Omega )$. We consider the
following two-phase inhomogeneous free boundary problem (f.b.p. in the
sequel). 
\begin{equation}
\begin{cases}
F(D^{2}u^{+})=f_{1} & \text{in }\Omega ^{+}(u):=\{u>0\} \\ 
F(D^{2}u^{-})=f_{2}\chi _{\{u<0\}} & \text{in }\Omega ^{-}(u)=\{u\leq 0\}^{o}
\\ 
u_{\nu }^{+}(x)=G(u_{\nu }^{-},x,\nu ) & \text{along }{\mathcal{F}}%
(u):=\partial \{u>0\}\cap \Omega .%
\end{cases}
\label{fbp}
\end{equation}%
Here $\nu =\nu (x)$ denotes the unit normal to the free boundary ${\mathcal{F%
}}={\mathcal{F}}(u)$ at the point $x$, pointing toward $\Omega ^{+}\left(
u\right) $, while the function $G(\beta ,x,\nu )$ is Lipschitz continuous,
strictly increasing in $\beta $, and 
\begin{equation}
\inf_{x\in \Omega ,|\nu |=1}G(0,x,\nu )>0.  \label{eq:G(0)}
\end{equation}%
Moreover, $u_{\nu }^{+}$ and $u_{\nu }^{-}$ denote the normal derivatives in
the inward direction to $\Omega ^{+}(u)$ and $\Omega ^{-}(u)$ respectively.

For any $u$ continuous in $\Omega $ we say that a point $x_{0}\in {\mathcal{F%
}}(u)$ is \emph{regular from the right} (resp. left) if there exists a ball $%
B\subset \Omega ^{+}(u)$ (resp. $B\subset \Omega ^{-}(u)$) such that $%
\overline{B}\cap {\mathcal{F}}(u)=x_{0}$. In both cases, we denote with $\nu
=\nu (x_{0})$ the unit normal to $\partial B$ at $x_{0}$, pointing toward $%
\Omega ^{+}(u)$.

\begin{definition}
\label{def:visc_sol_fbp} A \emph{viscosity} \emph{solution} of the free
boundary problem \eqref{fbp} is a continuous function $u$ which satisfies
the first two equality of (\ref{fbp}) in viscosity sense and such that the
free boundary condition is satisfied in the following viscosity sense:

(i) (supersolution condition) if $x_{0}\in \mathcal{F}(u)$ is regular from the
right with touching ball $B$, then, near $x_{0}$, 
\begin{equation*}
u^{+}(x)\geq \alpha \left\langle x-x_{0},\nu \right\rangle
^{+}+o(|x-x_{0}|)\qquad \text{in }B,\text{ with }\alpha \geq 0
\end{equation*}%
and 
\begin{equation*}
u^{-}(x)\leq \beta \left\langle x-x_{0},\nu \right\rangle
^{-}+o(|x-x_{0}|)\qquad \text{in }B^{c},\text{ with }\beta \geq 0,
\end{equation*}%
with equality along every non-tangential direction, and $\alpha \leq G(\beta
,x_{0},\nu (x_{0}));$

(ii) subsolution condition) if $x_{0}\in \mathcal{F}(u)$ is regular from the
left with touching ball $B$, then, near $x_{0}$, 
\begin{equation*}
u^{+}(x)\leq \alpha \left\langle x-x_{0},\nu \right\rangle
^{+}+o(|x-x_{0}|)\qquad \text{in }B^{c},\text{ with }\alpha \geq 0
\end{equation*}%
and 
\begin{equation*}
u^{-}(x)\geq \beta \left\langle x-x_{0},\nu \right\rangle
^{-}+o(|x-x_{0}|)\qquad \text{in }B,\text{ with }\beta \geq 0,
\end{equation*}%
with equality along every non-tangential direction, and 
\begin{equation*}
\alpha \geq G(\beta ,x_{0},\nu (x_{0})).
\end{equation*}
\end{definition}

The notion of viscosity solution can also be given in terms of test
functions (see [CS] for the equivalence). Given $u,\varphi \in C(\Omega )$,
we say that $\varphi $ touches $u$ by below (resp. above) at $x_{0}\in
\Omega ,$ if $u(x_{0})=\varphi (x_{0})$ and 
\begin{equation*}
u(x)\geq \varphi (x)\quad (\text{resp. $u(x)\leq \varphi (x)$})\quad \text{%
in a neighborhood $O$ of $x_{0}$.}
\end{equation*}%
Then $u\in C\left( \Omega \right) $ is a viscosity solution to (\ref{fbp})
if i) and ii) are replaced by

\smallskip

\textit{ii')} Let $x_{0}\in \mathcal{F}(u)$ and $v\in C^{2}(\overline{B^{+}(v)})\cap
C^{2}(\overline{B^{-}(v)})$ ($B=B_{\delta }(x_{0}),\delta $ small) with $%
\mathcal{F}(v)\in C^{2}$. If $v$ touches $u$ by below (resp. above) at $x_{0}$, then%
\begin{equation*}
v_{\nu }^{+}(x_{0})\leq G(v_{\nu }^{-},x_{0},\nu \left( x_{0}\right) ).\text{%
\ \ \ \ \ }(\text{resp. }\geq ).
\end{equation*}

\subsection{Existence of Lipschitz viscosity solutions and weak regularity
properties of the free boundary}


A solution to (\ref{fbp}) that involves partial differential equations that
are not in divergence form can be constructed via Perron's method, by taking
the infimum over the following class of \emph{admissible supersolutions} ${\ 
\mathcal{S}}$.

\begin{definition}
\label{class} A locally Lipschitz continuous function $w\in C(\overline{%
\Omega})$ is in the class ${\mathcal{S}}$ if

\begin{itemize}
\item[(a)] $w$ is a solution in viscosity sense to 
\begin{equation*}
\begin{cases}
F(D^2w^+)\le f_1 & \text{in }\Omega^+(w) \\ 
F(D^2w^-)\ge f_2\chi_{\{w<0\}} & \text{in }\Omega^-(w);%
\end{cases}%
\end{equation*}

\item[(b)] if $x_0\in {\mathcal{F}}(w)$ is regular from the left, with
touching ball $B$, then 
\begin{equation*}
w^+(x)\leq \alpha \left\langle x-x_0, \nu \right\rangle^+ + o(|x-x_0|) \qquad%
\text{in }B^c, \text{ with }\alpha\ge0
\end{equation*}
and 
\begin{equation*}
w^-(x)\geq \beta \left\langle x-x_0, \nu \right\rangle^- + o(|x-x_0|) \qquad%
\text{in }B, \text{ with }\beta\ge0,
\end{equation*}
with 
\begin{equation*}
\alpha \le G(\beta, x_0, \nu(x_0));
\end{equation*}

\item[(c)] if $x_0\in {\mathcal{F}}(w)$ is not regular from the left then 
\begin{equation*}
w(x)= o(|x-x_0|).
\end{equation*}
\end{itemize}
\end{definition}

The last ingredient one needs is that of minorant subsolution.

\begin{definition}
\label{minorante} A locally Lipschitz continuous function $\underline{u}\in
C(\overline{\Omega})$ is a \emph{strict minorant} if

\begin{itemize}
\item[(a)] $\underline{u}$ is a viscosity solution to 
\begin{equation*}
\begin{cases}
F(D^{2}\underline{u}^{+})\geq f_{1} & \text{in }\Omega ^{+}(\underline{u})
\\ 
F(D^{2}\underline{u}^{-})\leq f_{2}\chi _{\{\underline{u}<0\}} & \text{in }%
\Omega ^{-}(\underline{u});%
\end{cases}%
\end{equation*}

\item[(b)] every $x_{0}\in {\mathcal{F}}(\underline{u})$ is regular from the
right, with touching ball $B$, and near $x_{0}$ 
\begin{equation*}
\underline{u}^{+}(x)\geq \alpha \left\langle x-x_{0},\nu \right\rangle
^{+}+o(|x-x_{0}|)\qquad \text{in }B,\text{ with }\alpha >0,
\end{equation*}%
and 
\begin{equation*}
\underline{u}^{-}(x)\leq \beta \left\langle x-x_{0},\nu \right\rangle
^{-}+o(|x-x_{0}|)\qquad \text{in }B^{c},\text{ with }\beta \geq 0,
\end{equation*}%
with 
\begin{equation*}
\alpha >G(\beta ,x_{0},\nu (x_{0})).
\end{equation*}
\end{itemize}
\end{definition}

The following result holds (\cite{SVT}).

\begin{theorem}
\label{thm:main} Let $\,F$ be concave, homogeneous of degree one, and $g$ be
a continuous function on $\partial \Omega $. Assume that

\begin{itemize}
\item[(a)] there exists a strict minorant $\underline{u}$ with $\underline{u}%
=g$ on $\partial \Omega$ and

\item[(b)] the set $\{w\in {\mathcal{S}}:w\geq \underline{u},\ w=g\text{ on }%
\partial \Omega \}$ is not empty.

Let 
\begin{equation*}
u=\inf \{w:w\in {\mathcal{S}},w\geq \underline{u}\}.
\end{equation*}%
Then $u\in C\left( \overline{\Omega }\right) $, $u=g$ on $\partial \Omega $
and it is a (minimal) viscosity solution of (\ref{fbp}). Moreover ~$u$ is
locally Lipschitz in $\Omega $ with non degenerate positive part:%
\begin{equation*}
u^{+}\left( x\right) \geq \alpha \text{ dist}\left( x,F\left( u\right)
\right) \text{ \ \ \ }\left( \alpha >0\right) .
\end{equation*}
\end{itemize}
\end{theorem}

Once existence of a solution is established, we turn to the analysis of the
weak regularity properties of the free boundary.

The free boundary ${\mathcal{F}}(u)$ has finite $(n-1)$-dimensional
Hausdorff measure. More precisely, there exists a universal constant $%
r_{0}>0 $ such that for every $r<r_{0},$ for every $x_{0}\in {\mathcal{F}}%
(u),$ 
\begin{equation*}
\mathcal{H}^{n-1}({\mathcal{F}}(u)\cap B_{r}(x_{0}))\leq cr^{n-1}.
\end{equation*}%
Moreover, the reduced boundary ${\mathcal{F}}^{\ast }(u)$ of $\Omega ^{+}(u)$
has positive density in $\mathcal{H}^{n-1}-$ measure at any point of $F(u)$,
i.e. for $r<r_{0}$, $r_{0}$ universal 
\begin{equation*}
\mathcal{H}^{n-1}({\mathcal{F}}^{\ast }(u)\cap B_{r}(x))\geq cr^{n-1},
\end{equation*}%
for every $x\in {\mathcal{F}}(u).$ In particular 
\begin{equation*}
\mathcal{H}^{n-1}({\mathcal{F}}(u)\setminus {\mathcal{F}}^{\ast }(u))=0.
\end{equation*}%
Using the strong regularity results below we deduce the following result.

\begin{corollary}
${\mathcal{F}}(u)$ is a $C^{1,\gamma }$ surface in a neighborhood of ${%
\mathcal{H}}^{n-1}$ a.e. point $x_{0}\in {\mathcal{F}}(u).$
\end{corollary}

Existence of a continuous viscosity solution through a Perron method has
been established for linear operators in divergence form in \cite{C3}
(homogeneous case) and in \cite{DFS2} (inhomogeneous case), and for a class
of concave operators in \cite{W3}. 
 The presence of both a right hand side and the
nonlinearity of the governing equation presents several delicate points
which require new arguments, significantly when no sign condition is posed
on the right-hand-sides $f_{1}$ and $f_{2}$.

When $F=F\left( D^{2}u,Du\right) ,$ with $F$\ non concave in the Hessian
matrix and even for the linear case $F=$Tr$\left( A\left( x\right)
D^{2}u\right) +b\left( x\right) \cdot \nabla u,$ both existence and weak
regularity remain open questions.\smallskip

We conclude this section by mentioning a different approach to prove
existence for a specific class of free boundary problems in divergence form
that can be found for example in \cite{T}.

\section{Lipschitz continuity and global solutions}

\label{sec3}

The Lipschitz continuity of solutions to (\ref{fbp}) is a crucial ingredient
in the study of the regularity of the free boundary $F\left( u\right) $.
Indeed it provides compactness to carry on a blow-up analysis around a point 
$x_{0}\in F\left( u\right) ,$ reducing the problem to the classification of
global Lipschitz solution. For instance, if $F\left( Du\right) =-\Delta u$, $%
f_{1}=f_{2}=0$ it is possible to classify global solutions $U$ as either
purely \emph{two-plane} functions ($\nu $ a unit direction), 
\begin{equation*}
U(x)=\alpha \langle x-x_{0},\nu \rangle ^{+}-\beta \langle x-x_{0},\nu
\rangle ^{-}\text{ \ \ \ \ \ \ }\alpha =G\left( \beta \right) \text{, }\beta
>0
\end{equation*}%
or \emph{one-phase solutions}, in case we have 
\begin{equation*}
U^{-}\equiv 0.
\end{equation*}%
In particular, if $u$ is a solution with $0\in F(u)$, then via a blow-up
analysis and flatness results (see below) either $u^{-}(x)=o(|x|)$ or $F(u)$
is $C^{1,\gamma }$ in a neighborhood of $0$ that is, the only types of
singular points are the ones that occur in the one-phase setting when the
negative phase is identically zero.

On points of the reduced boundary of the minimal Perron solutions, one can
conclude that in the case when the blow up is a one-phase solution, then it
is in fact a one-plane solution $\alpha \langle x-x_{0},\nu \rangle $, hence 
$C^{1,\gamma }$ regularity of the reduced boundary follows from flatness
results below.

For viscosity solutions in general, when the governing equation is a fully
non linear operator more robust arguments are required.

\subsection{Lipschitz continuity}

In \cite{DS1} De Silva and Savin proved the following result under the
assumptions that $G$ behaves like $t$ for $t$ large.

\begin{theorem}
\label{DS1} Let $u$ be a viscosity solution to (\ref{fbp}) $f_{1}=f_{2}=0,$
and assume that $G\in C^{2}\left( [0,\infty \right) )$ and 
\begin{equation}
G^{\prime }\left( t\right) \rightarrow 1,\text{ \ \ }G^{^{\prime \prime
}}\left( t\right) =O\left( 1/t\right) \text{ as }t\rightarrow \infty .
\label{G}
\end{equation}%
Then 
\begin{equation*}
\left\Vert \nabla u\right\Vert _{L^{\infty }\left( B_{1/2}\right) }\leq
C\left( 1+\left\Vert u\right\Vert _{L^{\infty }\left( B_{1}\right) }\right)
\end{equation*}%
with $C=C\left( n,\lambda ,\Lambda ,G\right) $ universal.
\end{theorem}

The heuristic behind the proof is that ``big gradients" force the free
boundary condition to become a no-jump condition for $\nabla u$ and then
interior $C^{1,\alpha }$ estimates for fully nonlinear equations provide
gradient estimates.

The dependence on $G$ in the constant above is determined by the rate of
convergence in the limit \eqref{G}. In particular \eqref{G} can be relaxed
to $G^{\prime }\in \lbrack 1-\delta ,1+\delta ]$ for large values of $t$. If 
$F$ is homogeneous of degree 1, then it suffices to require $%
G(t)/t\rightarrow c_{0}$ as $t\rightarrow \infty $ for some constant $c_{0}$.

The method of the proof still works in the presence of a non-zero right hand
side.

In dimension $n=2$ these results can be improved significantly (see Section
4.3.)

\subsection{Classification of global solutions}

In \cite{DS2} De Silva and Savin proved the following Liouville theorem for
global Lipschitz solutions.

\begin{theorem}
\label{T2} Let $f_{1}=f_{2}=0$ and $u$ be a globally Lipschitz viscosity
solution to \eqref{fbp} in $\mathbb{R}^{n}.$ Assume that 
\begin{equation}
\text{ }F\text{ is concave (or convex) and homogeneous of degree 1.}
\label{F}
\end{equation}%
Then either $u$ is a two plane-solution 
\begin{equation}
u=\alpha \langle x-x_{0},\nu \rangle ^{+}-\beta \langle x-x_{0},\nu \rangle
^{-}\quad \quad \mbox{with}\quad \alpha ,\beta >0,\alpha =G(\beta ),
\label{2p}
\end{equation}%
or 
\begin{equation}
u^{-}\equiv 0,  \label{1p}
\end{equation}%
which means that $u$ solves the one-phase problem for $F$.
\end{theorem}

As in the case of the Laplacian, the main consequence of Theorem \ref{T2} is
that it reduces the question of the regularity of the free boundary for the
two-phase problem to the classification of global blow-up solutions to the
one-phase problem. In particular, if $u$ is a solution to $\eqref{fbp}$ with 
$0\in \Gamma (u)$, then either $u^{-}(x)=o(|x|)$ or $F(u)$ is $C^{1,\gamma }$
in a neighborhood of $0.$

Theorem \ref{T2} can be extended to more general operators $F(D^{2}u,\nabla
u,u)$ if appropriate assumptions are imposed on $F$. For example this result
holds when the problem is governed by quasilinear equations of the type 
\begin{equation}
\sum_{i,j=1}^{n}a_{ij}\left( \frac{\nabla u}{|\nabla u|}\right) \,\,u_{ij}=0,
\label{pL}
\end{equation}%
with uniformly elliptic coefficients $a_{ij}\in C^{1}(S^{n-1})$. In \cite%
{DS1} Lipschitz continuity of solutions to such a problem is also
established.

The proof requires somewhat involved and technical arguments. One of the
main steps consists in obtaining a weak Evans-Krylov type estimate for a
nonlinear transmission problem. One major difficulty is that $F(u)$ is not
known to be better than $C^{1,\alpha }$ even in the perturbative setting.

The idea of proof of Theorem \ref{T2} is to show a \textquotedblleft
reversed" improvement of flatness for the solution $u$, which means that if $%
u$ is sufficiently close to a two plane solution at a small scale then it
remains close to the same two-plane solution at all larger scales.

\subsection{Different operators}

When the two phases are governed by two different operators, few results are
known. In \cite{AT}, Amaral and Teixeira consider the case of divergence
form operators, with $f_{1},f_{2}\in L^{p},p>n/2$ and $A_{j}$ merely
measurable. They prove that local minimizers of the associated energy
functional have a universal H\"{o}lder modulus of continuity.

In \cite{CDS}, Caffarelli, De Silva and Savin obtained Lipschitz continuity
and classification of global Lipschitz solutions, that is Theorem \ref{DS1}
and Theorem \ref{T2}, for a two-phase problem driven by two different
operators with measurable coefficients under very general free boundary
conditions $u_{\nu }^{+}=G(u_{\nu }^{-},\nu ,x),$ in dimension $n=2$. We
remark that Theorem \ref{DS1} cannot hold in this generality in higher
dimensions. Indeed, say for $n=3$, it is not difficult to construct two
homogeneous functions of degree less than one, that solve two different
uniformly elliptic equations in complementary domains in $\mathbb{R}^{3}$
and satisfy the free boundary condition for a specific $G$.

\section{$C^{1,\alpha}$ regularity of the free boundary}

\label{sec4}

\subsection{The homogeneous case}

The regularity theory for the Laplace operator in the homogeneous case has
been developed by Caffarelli in the two seminal papers \cite{C1,C2}.

In particular the \textquotedblleft \emph{Lipschitz implies }$C^{1,\gamma }$%
\textquotedblright \emph{\ }part is contained in \cite{C1} while the \emph{%
flat implies Lipschitz }part is shown in \cite{C2}. The flatness condition
in \cite{C2} is stated in terms of \textquotedblleft $\varepsilon $-
monotonicity\textquotedblright\ along a cone of directions $\Gamma (\theta
_{0},e)$ of axis $e$ and large opening $\theta _{0}$. Precisely, a function $%
u$ is said to be $\varepsilon -$monotone ($\varepsilon >0$ small) along the
direction $\tau $ in the cone $\Gamma (\theta _{0},e)$ if for every $%
\varepsilon ^{\prime }\geq \varepsilon $, 
\begin{equation*}
u(x+\varepsilon ^{\prime }\tau )\leq u(x).
\end{equation*}%
Geometrically, the $\varepsilon $-monotonicity of $u$ can be interpreted as $%
\varepsilon $-closeness of $F(u)$ to the graph of a Lipschitz function.

In those papers Caffarelli set up a general strategy to attack the
regularity of the free boundary. Let us briefly describe the central idea of
the proof in \cite{C1}.

Starting from a Lipschitz graph, one shows that in a neighborhood of $%
F\left( u\right) $ the level sets of $u$ are still Lipschitz graphs, locally
in the same direction. Then one improves the Lipschitz constant (i.e. the
flatness) of the level sets of $u$ away from the free boundary. Here Harnack
inequality applied to directional derivatives of $u$ plays a major role.
Then the task is to carry this interior gain up to the free boundary. To
this aim, Caffarelli introduces a powerful method of continuity based on the
construction of a continuous family of deformations, constructed as the
supremum of a harmonic function over balls of variable radius (\emph{%
supconvolutions}). Finally, by rescaling and iterating the last two steps,
one obtains a geometric decay of the Lipschitz constant, which amounts to
the $C^{1,\gamma }$ regularity of $\mathcal{F}\left( u\right) $.

After 10 years Feldman \cite{F1} considered different anisotropic operators
with constant coefficients and extends to this case the results in \cite{C1}.

P. Y. Wang managed to extend the results in \cite{C1,C2} to a class of
concave fully non linear operators (see \cite{W1, W2}). One year later
Feldman in \cite{F2} considered a class of non concave fully non linear
operators of the type $F\left( D^{2}u,Du\right) $. He showed that Lipschitz
free boundaries are $C^{1,\alpha }$ thus extending to this case the results
in \cite{C1}.

The first papers dealing with variable coefficient operators are by Cerutti,
Ferrari, Salsa \cite{CFS1} and by Ferrari \cite{Fe} and Argiolas, Ferrari 
\cite{AF}. They considered respectively, linear elliptic operators in
non-divergence form and a rather general class of fully nonlinear operators $%
F\left( D^{2}u,Du,x\right) ,$ with H\"{o}lder continuity in $x$, including
Bellman's operators. One of the main difficulty in extending the theory to
variable coefficients operator is the fact that directional derivatives do
not satisfy any reasonable elliptic equation.

A refinement of the techniques in \cite{CFS1} leads to the following results
(see \cite{FeSa}), where the drift coefficient is merely bounded measurable,
with two different operators 
\begin{equation*}
F_{i}(D^{2}u,Du)=\emph{Tr}\left( A_{i}\left( x\right) D^{2}u\right)
+b_{i}\left( x\right) \cdot \nabla u\emph{,\ }i=1,2
\end{equation*}
governing the two phases.

\begin{theorem}
Let $.$ Let $u$ be a weak solution of our free boundary problem in $\Omega
=B_{1}=B_{1}\left( 0\right) ,$ the unit ball centered at the origin. Suppose 
$0\in F\left( u\right) $ and that

\begin{itemize}
\item[i)] $A_{i}\in C^{0,a}(B_{1})$\emph{, }$0<a\leq 1$\emph{, }$b_{i}\in
L^{\infty }\left( B_{1}\right) $\emph{.}

\item[ii)] $0<\alpha _{1}\leq \frac{u^{+}\left( x\right) }{dist\left(
x,F\left( u\right) \right) }\leq \alpha _{2}$\emph{.}

\item[iii)] $G\left( 0\right) >0.$\emph{\ }
\end{itemize}

There exist $0<\overline{\theta }<\pi /2\,$ and $\overline{\varepsilon }>0$
such that, if for $0<\varepsilon <\overline{\varepsilon },$ $\mathcal{F}\left(
u\right) $ is contained in an $\varepsilon -$neighborhood of a graph of a
Lipschitz function $x_{n}=g\left( x^{\prime }\right) $ with%
\begin{equation*}
Lip\left( g\right) \leq \tan \left( \frac{\pi }{2}-\overline{\theta }\right)
\end{equation*}%
then $\mathcal{F}\left( u\right) $ is a $C^{1,\gamma }-$graph in $B_{1/2}$.

The same conclusion holds if $B_1^{+}(u)=\{(x^{\prime
},x_{n}):x_{n}>f(x^{\prime })\}\cap B_{1}$\textsl{\ }where $f$\ is a
Lipschitz continuous function.
\end{theorem}

Condition $ii)$ expresses a linear behavior of $u^{+}$ at the free boundary
while being trapped in a neighborhood of two Lipschitz graph with small
Lipschitz constant is another way to express a flatness condition. Thus, 
\emph{flatness} plus \emph{linear behavior} of the positive part imply
smoothness.\medskip

This theorem has been extended to the same class of fully nonlinear
operators considered in \cite{Fe} by Argiolas and Ferrari in \cite{AF}. Note
that, in principle, all these results \emph{do not need the a priory
assumption of Lipschitz continuity of the solution}, that comes as a
consequence of the regularity of the free boundary.

\subsection{The nonhomogeneous case}

In this section we describe the strategy to investigate free boundary
problems with right hand side. Based on a Harnack type theorem and
linearization, this technique avoids the use of supconvolutions, that in
presence of distributed sources produces several complicacies. The method
can be very well adapted to nonhomogeneous two-phase problems to prove that
flat (see below) or Lipschitz free boundaries of (\ref{fbp}) are $%
C^{1,\gamma }.$ We assume that $f_{1}=f_{2}=f.$

We have (we will always assume that $0\in \mathcal{F}\left( u\right) $):

\begin{theorem}
\label{flatmain2} Let $u$ be Lipschitz viscosity solution to $(\ref{fbp})$
in $B_{1}.$ Assume that $f$ is bounded and continuous in $B_{1}^{+}(u)\cup
B_{1}^{-}(u).$ There exists a universal constant $\bar{\delta}>0$ such that,
if 
\begin{equation}
\{x_{n}\leq -\delta \}\subset B_{1}\cap \{u^{+}(x)=0\}\subset \{x_{n}\leq
\delta \},\ \ \ \ \ \ (\delta -\text{flatness})  \label{flat}
\end{equation}%
with $0\leq \delta \leq \bar{\delta},$ then $F(u)$ is $C^{1,\gamma }$ in $%
B_{1/2}$.
\end{theorem}

Condition (\ref{flat}) expresses that the zero set of $u^{+}$ is trapped
between two parallel hyperplanes at $\delta -$distance from each other for a
small $\delta $, which is a kind of flatness ($\delta -$flatness)$.$ While
this looks like a somewhat strong assumption, it is indeed a natural one
since it is satisfied for example by rescaling a solution around a point of
the free boundary where there is a normal in some weak sense (\emph{regular
points}), for instance in the measure theoretical one. Thus, for the minimal
Perron solution $H^{n-1}-$a.e. points on $\mathcal{F}\left( u\right) $ are
of this kind. Moreover, starting from a Lipschitz free boundary, $H^{n-1}-$%
a.e. points on $\mathcal{F}\left( u\right) $ are regular, by Rademacher\
Theorem. \smallskip

When $F$ is (positively) homogeneous of degree one (or when $F_{r}(M)=\frac{1%
}{r}F\left( rM\right) $ has a limit $F^{\ast }(M),$ as $%
r\rightarrow 0,$ which is always homogeneous of degree one), we also have:

\begin{theorem}
\label{Lipmainvar} (Lipschitz implies $C^{1,\gamma }$)\label{LIP} Let $u$ be
a Lipschitz viscosity solution to $(\ref{fbp})$ in $B_{1}.$ Assume that $f$
is bounded and continuous in $B_{1}^{+}(u)\cup B_{1}^{-}(u).$ If $\mathcal{F}%
(u)$ is a Lipschitz graph in a neighborhood of $0,$ then $\mathcal{F}(u)$ is 
$C^{1,\gamma }$in a (smaller) neighborhood of $0$.
\end{theorem}

Theorem \ref{LIP} follows from  Theorem \ref{flatmain2} via a blow-up argument and a Liouville type result for global viscosity solutions to a two-phase free boundary problem, for the details see the original paper \cite{DFS1}.\smallskip

As already pointed out the proof of Theorem \ref{flatmain2} follows the
strategy developed in \cite{D}. The main difficulty in the analysis in this
two-phase problem comes from the case when $u^{-}$ is \emph{degenerate},
that is very close to zero without being identically zero. In this case the
flatness assumption does not guarantee closeness of $u$ to an
\textquotedblleft optimal\textquotedblright\ (two-plane or one-plane)
configuration. Thus one needs to work only with the positive phase $u^{+}$
to balance the situation in which $u^{+}$ highly predominates over $u^{-}$
and the case in which $u^{-}$ is not too small with respect to $u^{+}.$

In particular, the proof is based on a recursive improvement of flatness,
obtained via a compactness argument, provided by a geometric type Harnack
inequality, which linearizes the problem into a limiting one. The limiting
problem turns out to be a transmission problem in the nondegenerate case and
a Neumann problem in the other case. The information to set up the iteration
towards regularity is precisely stored in the analysis of this problem.

Let us heuristically show how the free boundary condition 
\begin{equation*}
\left\vert \nabla u^{+}\right\vert =G\left( \left\vert \nabla
u^{-}\right\vert \right)
\end{equation*}%
linearizes in the \emph{nondegenerate }case.\emph{\ }Let $U_{\beta
}\left( t\right) =\alpha t^{+}-\beta t^{-},$ $\alpha =G(\beta ),$ it is possible to show that for any $\varepsilon>0,$ if $\delta>0,$ depending on $\varepsilon$ is small enough, the condition \eqref{flat} implies the following one 
\begin{equation}
U_{\beta }(x_{n}-\varepsilon )\leq u(x)\leq U_{\beta }(x_{n}+\varepsilon
)\quad \text{in $B_{1},$}  \label{ff1}
\end{equation}%
with $0<\beta \leq L$, and $L=$Lip$\left( u\right) $. After rescaling we may
assume that 
\begin{equation*}
\left\vert f\right\vert \leq \varepsilon ^{2}\min \left\{ \alpha ,\beta
\right\}
\end{equation*}

This suggests the renormalization 
\begin{equation*}
\tilde{u}_{\varepsilon }(x)=%
\begin{cases}
\dfrac{u(x)-\alpha x_{n}}{\alpha \varepsilon },\quad x\in B_{1}^{+}(u)\cup 
\mathcal{F}(u) \\ 
\  \\ 
\dfrac{u(x)-\beta x_{n}}{\beta \varepsilon },\quad x\in B_{1}^{-}(u)%
\end{cases}%
\end{equation*}%
or%
\begin{equation}
u\left( x\right) =%
\begin{cases}
\alpha x_{n}+\varepsilon \alpha \tilde{u}_{\varepsilon }(x),\quad x\in
B_{1}^{+}(u)\cup \mathcal{F}(u) \\ 
\  \\ 
\beta x_{n}+\varepsilon \beta \tilde{u}_{\varepsilon }(x),\quad x\in
B_{1}^{-}(u).%
\end{cases}
\label{pert}
\end{equation}%
Recalling \eqref{operatorname}, we have
\begin{equation*}
F_{\alpha\varepsilon}(D^2\tilde{u}_{\varepsilon })=\frac{f}{\alpha \varepsilon 
}\sim \varepsilon \text{ \ \ \ \ in }B_{1}^{+}\left( u\right).
\end{equation*}
and
\begin{equation*}
	F_{\beta\varepsilon}(D^2\tilde{u}_{\varepsilon })=\frac{f}{\beta \varepsilon 
	}\sim \varepsilon \text{ \ \ \ \ in }B_{1}^{-}\left( u\right).
\end{equation*}

On $\mathcal{F}\left( u\right) ,$%
\begin{equation*}
\left\vert \nabla u^{+}\right\vert =\alpha \left\vert e_{n}+\varepsilon
\nabla \tilde{u}_{\varepsilon }(x)\right\vert \sim \alpha \left(
1+\varepsilon \left( \tilde{u}_{\varepsilon }\right) _{x_{n}}+\varepsilon
^{2}\left\vert \nabla \tilde{u}_{\varepsilon }\right\vert ^{2}\right)
\end{equation*}%
and 
\begin{eqnarray*}
G\left( \left\vert \nabla u^{-}\right\vert \right) &=&G\left( \left\vert
\beta e_{n}+\varepsilon \beta\nabla \tilde{u}_{\varepsilon }\right\vert \right) \sim
G\left(\beta\left( 1+\varepsilon \left( \tilde{u}_{\varepsilon }\right)
_{x_{n}}+\varepsilon ^{2}\left\vert \nabla \tilde{u}_{\varepsilon
}\right\vert ^{2}\right) \right) \\
&\sim &G(\beta )+\varepsilon G^{\prime }\left( \beta \right) \left( \beta
\left( \tilde{u}_{\varepsilon }\right) _{x_{n}}+\varepsilon \beta \left\vert
\nabla \tilde{u}_{\varepsilon }\right\vert ^{2}\right) +\varepsilon ^{2}.
\end{eqnarray*}%
Letting $\varepsilon \rightarrow 0$, we get formally for \textquotedblleft
the limit\textquotedblright\ $\tilde{u}$ the following problem: 
\begin{equation}
F^{\pm }\left( D^{2}\tilde{u}\right) =0\text{, \ in }B_{1/2}^{\pm }\cap
\left\{ x_{n}\neq 0\right\}  \label{tr}
\end{equation}%
and the transmission condition (\emph{linearization of the free boundary
condition})%
\begin{equation}
\alpha \left( \tilde{u}_{x_{n}}\right) ^{+}-\beta G^{\prime }\left( \beta\right)
\left( \tilde{u}_{x_{n}}\right) ^{-}=0\text{ \ on }B_{1/2}\cap \left\{
x_{n}=0\right\}  \label{trcon}
\end{equation}%
where $F^{+}\left( M\right) ,F^{-}\left( M\right) $ are limits (of
sequences) of operators of the form 
\begin{equation*}
F_{\alpha \varepsilon }\left(M\right) \text{
\ and } F_{\beta \varepsilon }\left(M\right) ,%
\text{\ }
\end{equation*}%
while $\left( \tilde{u}_{x_{n}}\right) ^{+}$ and $\left( \tilde{u}%
_{x_{n}}\right) ^{-}$ denote the $e_{n}-$derivatives of $\tilde{u}$
restricted to $\{x_{n}>0\}$ and $\{x_{n}<0\}$, respectively.

Thus, at least formally, we have found an asymptotic problem for the limits
of the renormalizations $\tilde{u}_{\varepsilon }$. The crucial information
we were mentioning before is contained in the following regularity result
that we state for the fully nonlinear case (see \cite{DFS5} when distributed
sources are present). Consider the transmission problem, ($\tilde{a}\neq 0$) 
\begin{equation}
\begin{cases}
F^{+}\left( D^{2}\tilde{u}\right) =0 & \text{in $B_{1}^{+}\cap
\{x_{n}\neq 0\}$}, \\ 
F^{-}\left( D^{2}\tilde{u}\right)=0 & \text{in $B_{1}^{-}\cap
\{x_{n}\neq 0\}$}, \\ 
\tilde{a}\left( \tilde{u}_{x_{n}}\right) ^{+}-\tilde{b}\left( \tilde{u}_{x_{n}}\right) ^{-}=0 & \text{on $%
B_{1}\cap \{x_{n}=0\}$}.%
\end{cases}
\label{Neu}
\end{equation}

\begin{theorem}
\label{linearreg}Let $\tilde{u}$ be a viscosity solution to $\eqref{Neu}$ in 
$B_{1}$ such that $\Vert \tilde{u}\Vert _{\infty }\leq 1.$ Then $\tilde{u}%
\in C^{1,\gamma}\left( \bar{B}_{1}^{\pm }\right) $ and in particular, there
exists a universal constant $\bar{C}$ such that 
\begin{equation}
|\tilde{u}(x)-\tilde{u}(0)-(\nabla _{x^{\prime }}\tilde{u}(0)\cdot x^{\prime
}+\tilde{p}x_{n}^{+}-\tilde{q}x_{n}^{-})|\leq \bar{C}r^{1+\gamma },\quad 
\text{in $B_{r}$}  \label{up}
\end{equation}%
for all $r\leq 1/2,$ with $\tilde{a}\tilde{p}-\tilde{b}\tilde{q}=0.$
\end{theorem}

Transferring the estimate (\ref{up}) to $\tilde{u}_{\varepsilon }$ and then
reading it in terms of flatness for $u$ through formulas (\ref{pert}), one
deduce that if $0<r\leq r_{0}$ for $r_{0}$ universal, and $0<\varepsilon
\leq \varepsilon _{0}$ for some $\varepsilon _{0}$ depending on $r$, 
\begin{equation}
U_{\beta^{\prime }}(x\cdot \nu _{1}-r\frac{\varepsilon }{2})\leq u(x)\leq
U_{\beta^{\prime }}(x\cdot \nu _{1}+r\frac{\varepsilon }{2})\quad \text{in }%
B_{r},  \label{upup}
\end{equation}
with $|\nu _{1}|=1,$ $|\nu _{1}-e_{n}|\leq \tilde{C}\varepsilon $ , and $%
|\beta-\beta^{\prime }|\leq \tilde{C}\beta\varepsilon $ for a universal
constant $\tilde{C}.$

Rescaling and iterating (\ref{upup}) one deduces uniform pointwise $%
C^{1,\gamma }$ estimates in \ neighborhood of the origin.

\section{Higher regularity}

\label{sec5}

In view of the results in Section 5, under suitable flatness assumptions,
the free boundary $\mathcal{F}(u)$ is locally $C^{1,\gamma }$ and the same
conclusion holds if $\mathcal{F}\left( u\right) $ is a graph of a Lipschitz
function. Therefore $u$ is a classical solution, i.e. the free boundary
condition is satisfied in a pointwise sense. Concerning higher regularity,
the result one aims to prove should be the following one.
\smallskip

"\textbf{Theorem}". Assume 
 some reasonable starting regularity on $\mathcal{F}\left( u\right) $ and let $k$ be a nonnegative integer. Suppose that 
$f_{i}\in C^{k,\gamma }\left( B_{1}\right) $, $i=1,2,$
and $G$ is $C^{2+k}.$ Then $\mathcal{F}\left( u\right) \cap B_{1/2}
$ is $C^{k+2,\gamma ^{\ast }}.$ Moreover, if $f_{i}$, $i=1,2,$ are
$C^{\infty }$ or real analytic in $B_{1}$, then $
\mathcal{F}\left( u\right) \cap B_{1/2}$ is $C^{\infty }$ or real
analytic, respectively.\smallskip 

The state of the art is the following. In the seminal paper \cite{KNS}, the
authors used a zero order hodograph transformation and a suitable reflection
map, to locally reduce a two-phase problem to an elliptic, coercive
nonlinear system of equations. The existing literature on the regularity of
solutions to nonlinear systems developed in \cite{ADN, Morrey} can be
applied as long as the solution $u$ is $C^{2,\alpha }$ for some $\alpha >0$
up to the free boundary (from either side).

As noted in the recent work \cite{KL}, in the case when the governing
equation in \eqref{fbp} is linear and in divergence form, the initial
assumption to obtain the above theorem is actually $u\in C^{1,\alpha }$. It
is not evident that the case of linear nondivergence uniformly elliptic
equations can be treated in a similar manner. On the other hand, the case
when the leading operator is a fully nonlinear operator definitely requires
the solution to have H\"{o}lder second derivatives (from both sides).

In \cite{DFS3} we consider the case of nondivergence linear operators with
the purpose to develop a general strategy that would apply to a larger class
of problem, to include also fully nonlinear operators. The application to
the latter case would be rather straightforward once $C^{2,\alpha }$
estimates for the limiting problem (\ref{Neu}) were available. However for
fully nonlinear operators, this remains an open problem. The main result in 
\cite{DFS3} is the following.

\begin{theorem}
\label{flatmain1} Let $F$ be a linear uniformly elliptic operator in  nondivergence form, with $C^{0,\gamma}$ coefficients and let $u$ be
a (Lipschitz) viscosity solution to \eqref{fbp} in $B_{1}.$ Assume that either $\mathcal{F}(u)$ is locally the graph of a Lipschitz continuous function or there exists a
universal constant $\eta>0$ such that, 
\begin{equation}
\{x_{n}\leq -\eta \}\subset B_{1}\cap \{u^{+}(x)=0\}\subset \{x_{n}\leq \eta
\},
\end{equation} 
then $\mathcal{F}(u)$ is $C^{2,\gamma ^{\ast }}$ in $B_{1/2}$ for a small $\gamma
^{\ast }$ universal, with the $C^{2,\gamma ^{\ast }}$ norm bounded by a
universal constant.
\end{theorem}
Hence, the initial proposal, stated in the theoretical expected "Theorem", may be considered achieved by the following corollary contained in \cite{DFS3}.
\begin{corollary}\label{corcor} Let $u$ be a viscosity solution of \eqref{fbp} when $F:=\mbox{Tr}(A(x)D^2u(x)),$  and $A$ is a smooth uniformly elliptic matrix of coefficients.
  Assume that $\mathcal{F}\left( u\right) $ is Lipschitz 
and let $k$ be a nonnegative integer. Suppose that 
$f_{i}\in C^{k,\gamma }\left( B_{1}\right) $, $i=1,2,$
and $G$ is $C^{2+k}.$ Then $\mathcal{F}\left( u\right) \cap B_{1/2}
$ is $C^{k+2,\gamma ^{\ast }}.$ Moreover, if $f_{i}$, $i=1,2,$ are
$C^{\infty }$ or real analytic in $B_{1}$, then $
\mathcal{F}\left( u\right) \cap B_{1/2}$ is $C^{\infty }$ or real
analytic, respectively.
\end{corollary}
 
Other related higher regularity results can be found in \cite{E, K}. The
overall strategy for the proof of Theorem \ref{flatmain1} follows the ideas
described in the previous section. However, reaching the $C^{2,\gamma }$
regularity requires a much more involved process because of the possible
degeneracy of the negative part. Indeed this causes a delicate interplay
between the two phases. Ultimately the main source of difficulties is due to
the presence of a forcing term of general sign in the negative phase. Indeed, if $f_{2}\geq 0$, the Hopf maximum principle would imply
nondegeneracy (also) on the negative side, making the two-phases of
comparable size and considerably simplifying the final iteration procedure.
It is worth noticing, that however even in this easier scenario (and in
particular in the homogeneous case), if one wants to attain uniform
estimates with universal constants, then one must employ the more involved
methods developed in \cite{DFS3} for the degenerate case.

\section{Nonlinear operators with non-standard growth}

\label{sec6}

In the second part of this survey we present some recent results concerning a
one-phase free boundary problem governed by the $p(x)$-Laplacian. More
precisely, let 
\begin{equation*}
\Delta _{p(x)}u=\mbox{div}(|\nabla u|^{p(x)-2}\nabla u),
\end{equation*}%
where $p$ is measurable and $1<p(x)<+\infty $
a.e. in a domain $\Omega $. We recall that if $p:\Omega
\rightarrow \lbrack 1,\infty )$ is a measurable bounded function,
the variable exponent Lebesgue space $L^{p(\cdot )}(\Omega )$ is defined as
the set of all measurable functions $u:\Omega \rightarrow \mathbb{R}$ for
which the modular $\varrho _{p(\cdot )}(u)=\int_{\Omega }|u(x)|^{p(x)}\,dx$
is finite. The Luxemburg norm on this space is defined by 
\begin{equation*}
\Vert u\Vert _{L^{p(\cdot )}(\Omega )}=\Vert u\Vert _{p(\cdot )}=\inf
\{\lambda >0:\varrho _{p(\cdot )}(u/\lambda )\leq 1\}.
\end{equation*}%
This norm makes $L^{p(\cdot )}(\Omega )$ a Banach space. In addition,
setting $p_{\max }=\mathrm{esssup}\,p(x)$ and $p_{\min }=\mathrm{essinf}%
\,p(x),$  the following relation holds between $\varrho _{p(\cdot )}(u)$
and $\Vert u\Vert _{L^{p(\cdot )}}$: 
\begin{equation*}
\begin{split}
& \min \Big\{\Big(\int_{\Omega }|u|^{p(x)}\,dx\Big)^{1/{p_{\min }}},\Big(%
\int_{\Omega }|u|^{p(x)}\,dx\Big)^{1/{p_{\max }}}\Big\}\leq \Vert u\Vert
_{L^{p(\cdot )}(\Omega )} \\
& \leq \max \Big\{\Big(\int_{\Omega }|u|^{p(x)}\,dx\Big)^{1/{p_{\min }}},%
\Big(\int_{\Omega }|u|^{p(x)}\,dx\Big)^{1/{p_{\max }}}\Big\}.
\end{split}%
\end{equation*}

Moreover, the dual of $L^{p(\cdot)}(\Omega)$ is $L^{p^{\prime
}(\cdot)}(\Omega)$ with $\frac{1}{p(x)}+\frac{1}{p^{\prime }(x)}=1$.

$W^{1,p(\cdot)}(\Omega)$ denotes the space of measurable functions $u$ such
that $u$ and the distributional derivative $\nabla u$ are in $%
L^{p(\cdot)}(\Omega)$. The norm

\begin{equation*}
\Vert u\Vert _{1,p(\cdot )}:=\Vert u\Vert _{p(\cdot )}+\Vert |\nabla u|\Vert
_{p(\cdot )}
\end{equation*}%
makes $W^{1,p(\cdot )}(\Omega )$ a Banach space. The space $W_{0}^{1,p(\cdot
)}(\Omega )$ is defined as the closure of the $C_{0}^{\infty }(\Omega )$ in $%
W^{1,p(\cdot )}(\Omega )$.

Assume that $\Omega \subset {\mathbb{R}}^{n}$ is a bounded domain, 
$p\in C^{1}(\Omega )$, $f\in C(\Omega )\cap L^{\infty }(\Omega )$ and $g\in
C^{0,\beta }(\Omega ),$ $g\geq 0.$ Then our problem reads 
\begin{equation}
\left\{ 
\begin{array}{ll}
\Delta _{p(x)}u=f, & \hbox{in $\Omega^+(u):= \{x \in \Omega : u(x)>0\}$}, \\ 
\  &  \\ 
|\nabla u|=g, & \hbox{on $F(u):= \partial \Omega^+(u) \cap
\Omega.$}%
\end{array}%
\right.   \label{fb}
\end{equation}

In addition we assume that 
\begin{equation}
\nabla p\in L^{\infty }(\Omega )  \label{p-lip}
\end{equation}%
and that 
\begin{equation}
1<p_{\min }\leq p(x)\leq p_{\max }<\infty .  \label{exponentsize}
\end{equation}

This kind of problems arise naturally from limits of a singular perturbation
problem with forcing term as in \cite{LW1}, where the authors analyze
solutions to \eqref{fb}, arising in the study of flame propagation with
nonlocal and electromagnetic effects. On the other hand, \eqref{fb} appears
by minimizing the energy functional 
\begin{equation}
\mathcal{E}(v)=\int_{\Omega }\left( \frac{|\nabla v|^{p(x)}}{p(x)}%
+Q^{2}(x)\chi _{\{v>0\}}+f(x)v\right) dx  \label{AC-energy}
\end{equation}
studied in \cite{LW3}, see \cite{DF} as well. We refer also to \cite{LW4}, where \eqref{fb}
appears in the study of an optimal design problem.

In general, partial differential equations with non-standard growth have
been receiving a lot of attention and the $p(x)$-Laplacian is a model case
in this class. A list of applications of this type of operators includes the
modelling of non-Newtonian fluids, for instance, electrorheological \cite{R}
or thermorheological fluids \cite{AR}. Also non-linear elasticity \cite{Z1},
image reconstruction \cite{AMS,CLR} and the modelling of electric conductors 
\cite{Z2}, to cite some of them.

\bigskip

We are interested in the regularity of the free boundary for
viscosity solutions of (\ref{fb}), following the strategy introduced in 
\cite{D}. The same technique was applied to the $p$-Laplace operator 
($p(x)\equiv p$ in \eqref{fb}), with $p\geq 2$, in \cite {LR}. 

Problem \eqref{fb} was originally studied in the linear homogeneous case
in \cite{AC}, associated to \eqref{AC-energy}. These techniques were
generalized to the linear case with $f\not\equiv 0$ in \cite{GS, Le} and, in
the homogeneous case, to a quasilinear uniformly elliptic situation \cite%
{ACF}, to the $p$-Laplacian \cite{DP}, to an Orlicz setting \cite{MW} and to
the $p(x)$-Laplacian with $p(x)\geq 2$ \cite{FMW}. For the case \eqref{fb}
with $1<p(x)<\infty $ and $f\not\equiv 0$ we refer to \cite{LW2}.


\bigskip


We will  present results on problem  \eqref{fb} obtained in \cite{FL}. In order to do so, we start with some definitions and preliminaries. 

Recall that we have assumed that $1<p_{\min }\leq p(x)\leq p_{\max }<\infty $, with $p(x)$
Lipschitz continuous in $\Omega$, and $f\in L^{\infty }(\Omega )$. We say that $u$
is a \emph{weak} solution to $\Delta _{p(x)}u=f$ in $\Omega $ if $u\in
W^{1,p(\cdot )}(\Omega )$ and, for every $\varphi \in C_{0}^{\infty }(\Omega
)$, there holds that 
\begin{equation}
-\int_{\Omega }|\nabla u(x)|^{p(x)-2}\nabla u\cdot \nabla \varphi
\,dx=\int_{\Omega }\varphi \,f(x)\,dx.  \label{defnweak}
\end{equation}%

The notion of viscosity solution to \eqref{fb} is given as
follows.

\begin{definition}
\label{defnhsol1} Let $u$ be a continuous nonnegative function in $\Omega$.
We say that $u$ is a viscosity solution to (\ref{fb}) in $\Omega$, if the
following conditions are satisfied:

\begin{enumerate}
\item $\Delta_{p(x)} u = f$ in $\Omega^+(u)$ in the weak sense (i.e., in the sense of \eqref{defnweak}).

\item For every $\varphi \in C(\Omega )$, $\varphi \in C^{2}(\overline{%
\Omega ^{+}(\varphi )})$. If $\varphi ^{+}$ touches $u$ from below (resp.
above) at $x_{0}\in F(u)$ and $\nabla \varphi (x_{0})\not=0$, then 
\begin{equation*}
|\nabla \varphi (x_{0})|\leq g(x_{0})\quad (\text{resp. $\geq g(x_{0})$}).
\end{equation*}
\end{enumerate}
\end{definition}

\smallskip

We point out the difference between the notion of viscosity solution for free boundary problem \eqref{fb} in Definition \ref{defnhsol1}  and the one for  free boundary problem \eqref{fbp} in Definition \ref{def:visc_sol_fbp} (Section \ref{sec2}). Namely, equation $\Delta_{p(x)} u = f$ is satisfied here  in weak the sense of \eqref{defnweak}. However, the program followed in \cite{FL} required to deal with
viscosity solutions of this equation --in the sense of \cite{CIL}-- as well. 

The equivalence between weak and viscosity solutions of $\Delta _{p(x)}u=f$
was proved in \cite{JJ,JLM, MO} in the case of the $p$-Laplacian (i.e., for $%
p(x)\equiv p$) and in \cite{JLP} in the case of the homogeneous $p(x)$%
-Laplacian (i.e., for $f\equiv 0$). For the inhomogeneous $p(x)$-Laplacian the following result was proven in \cite{FL}:  

\begin{theorem}
\label{weak-is-visc copy} Assume that $f\in C(\Omega )\cap L^{\infty
}(\Omega ),$ $p\in C^{1}(\Omega )$ with $1<p_{\min }\leq p(x)\leq p_{\max
}<\infty $ and $\nabla p\in L^{\infty }(\Omega )$. Let $u\in
W^{1,p(\cdot )}(\Omega )\cap C(\Omega )$ be a weak solution to $\Delta
_{p(x)}u=f$ in $\Omega $. Then $u$ is a viscosity solution to $\Delta
_{p(x)}u=f$ in $\Omega .$
\end{theorem}

As a consequence, we have: 

\begin{proposition}
\label{defnhsol2} Let $u$ be a viscosity solution to (\ref{fb}) in $\Omega$. 
Then the following conditions are satisfied:

\begin{enumerate}
\item $\Delta _{p(x)}u=f$ in $\Omega ^{+}(u)$ in the viscosity sense, that
is:

\begin{itemize}
\item[(ia)] for every $\varphi \in C^{2}(\Omega ^{+}(u))$ and for every $%
x_{0}\in \Omega ^{+}(u),$ if $\varphi $ touches $u$ from above at $x_{0}$
and $\nabla \varphi (x_{0})\not=0,$ then $\Delta _{p(x_{0})}\varphi
(x_{0})\geq f(x_{0}),$ that is, $u$ is a viscosity subsolution;

\item[(ib)] for every $\varphi \in C^{2}(\Omega ^{+}(u))$ and for every $%
x_{0}\in \Omega ^{+}(u),$ if $\varphi $ touches $u$ from below at $x_{0}$
and $\nabla \varphi (x_{0})\not=0,$ then $\Delta _{p(x_{0})}\varphi
(x_{0})\leq f(x_{0}),$ that is, $u$ is a viscosity supersolution.
\end{itemize}

\item For every $\varphi \in C(\Omega )$, $\varphi \in C^{2}(\overline{%
\Omega ^{+}(\varphi )})$. If $\varphi ^{+}$ touches $u$ from below (resp.
above) at $x_{0}\in F(u)$ and $\nabla \varphi (x_{0})\not=0$, then 
\begin{equation*}
|\nabla \varphi (x_{0})|\leq g(x_{0})\quad \mbox{(resp.}\geq g(x_{0})\mbox{)}%
.
\end{equation*}
\end{enumerate}
\end{proposition}

We now describe the main result in \cite{FL}. Namely, flat free
boundaries of viscosity solutions to \eqref{fb} are $C^{1,\alpha }.$ In the forthcoming work \cite{FL2} it is shown that Lipschitz free
boundaries of viscosity solutions are $C^{1,\alpha}.$  Concerning
this last paper, we just say that, remarkably, the regularity step from
Lipschitz to $C^{1,\alpha }$ presents new and interesting aspects,
requiring deeper techniques. 

\smallskip

Precisely, the main result in \cite{FL} is the following 

\begin{theorem}[Flatness implies $C^{1,\protect\alpha}$]
\label{flatmain1fl} Let $u$ be a viscosity solution to \eqref{fb} in $B_1$.
Assume that $0\in F(u),$ $g(0)=1$ and $p(0)=p_0.$ There exists a universal
constant $\bar{\varepsilon}>0$ such that, if the graph of $u$ is $\bar{%
\varepsilon}-$flat in $B_1,$ in the direction $e_n,$ that is 
\begin{equation}  \label{cflat}
(x_n-\bar{\varepsilon})^+\leq u(x)\leq (x_n+\bar{\varepsilon})^+, \quad x\in
B_1,
\end{equation}
and 
\begin{equation}  \label{pflat}
\|\nabla p\|_{L^{\infty}(B_1)}\leq \bar{\varepsilon},\quad
\|f\|_{L^{\infty}(B_1)}\leq \bar{\varepsilon}, \quad
[g]_{C^{0,\beta}(B_1)}\leq \bar{\varepsilon},
\end{equation}
then $F(u)$ is $C^{1,\alpha}$ in $B_{1/2}$.
\end{theorem}

\medskip

In Theorem \ref{flatmain1fl} the constants $\bar{\varepsilon}$ and $\alpha $
depend only on $p_{\min }$, $p_{\max }$ and $n$. 

\bigskip 

We  remark that Theorem \ref{flatmain1fl} is crucial in the
companion paper \cite{FL2} to prove that Lipschitz free boundaries of
viscosity solutions of \eqref{fb} are $C^{1,\alpha }$.

The proof of Theorem \ref{flatmain1fl} is based on an improvement
of flatness, obtained via a compactness argument (provided by a Harnack type
inequality) which linearizes the problem into a limiting one. We point out that the
development of new tools was necessary in order to implement this strategy for the
inhomogeneous $p(x)$-Laplace operator. In fact, the $p(x)$%
-Laplacian is a nonlinear operator that appears naturally in divergence form
from minimization problems, i.e., in the form $\mathrm{div}A(x,\nabla u)=f(x)
$, with 
\begin{equation*}
\lambda |\eta |^{p(x)-2}|\xi |^{2}\leq \sum_{i,j=1}^{n}\frac{\partial A_{i}}{%
\partial \eta _{j}}(x,\eta )\xi _{i}\xi _{j}\leq \Lambda |\eta
|^{p(x)-2}|\xi |^{2},\quad \xi \in \mathbb{R}^{n}.
\end{equation*}%
This operator is singular in the regions where $1<p(x)<2$ and degenerate in
the ones where $p(x)>2$.

Some important results for this type of operators  are available in the literature only for weak solutions (in the sense of \eqref{defnweak}). These results are Harnack inequality, \cite{Wo}, and $C^{1,\alpha}$ estimates, \cite{Fan}, \cite{FanZ}. This is the motivation for the choice of Definition \ref{defnhsol1}.
However, the proof of Theorem \ref{flatmain1fl} relies on solutions of inhomogeneous $p(x)$-Laplace equation in the viscosity sense. Then, Theorem \ref{weak-is-visc copy} becomes essential.


On the other hand, the
nondivergence nature of the viscosity solutions requires the construction of
suitable barriers. It turns out that barriers of the type $w(x)=c_{1}|x-x_{0}|^{-%
\gamma }-c_{2}$ are appropriate for the inhomogeneous $p(x)$-Laplace operator, although the proof is delicate,  due to
the nonlinear singular/degenerate nature of the operator, its $x$
dependence and the presence of the logarithmic term that appears in the nondivergence
form of the operator. Precisely, the following key result was proved in \cite{FL}: 

\begin{lemma}
Let $x_{0}\in B_{1}$ and $0<\bar{r}_{1}<\bar{r}_{2}\leq 1$. Assume that $%
1<p_{\min }\leq p(x)\leq p_{\max }<\infty $ and $\Vert \nabla p\Vert
_{L^{\infty }}\leq \varepsilon ^{1+\theta }$, for some $0<\theta \leq 1$.
Let $c_{0},c_{1},c_{2}$ be positive constants and $c_{3}\in \mathbb{R%
}$.
There exist positive constants $\gamma\ge 1$, $\bar{c}$, $\varepsilon_0$ and 
$\varepsilon_1$ such that the functions 
\begin{equation*}
w(x)=c_1|x-x_0|^{-\gamma}-c_2,
\end{equation*}
\begin{equation*}
v(x)=q(x)+\frac{c_0}{2}\varepsilon (w(x)-1),\quad q(x)=x_n+c_3
\end{equation*}
satisfy, for $\bar{r}_1\le |x-x_0|\le \bar{r}_2$, 
\begin{equation}  \label{eq-w}
\Delta_{p(x)}w\geq \bar{c}, \quad \text{for }\, 0<\varepsilon\le
\varepsilon_0,
\end{equation}
\begin{equation}  \label{eq-v}
\frac{1}{2}\le |\nabla v|\le 2, \qquad \Delta_{p(x)}v > \varepsilon^2, \quad 
\text{for }\, 0<\varepsilon\le \varepsilon_1.
\end{equation}
Here $\gamma=\gamma(n,p_{\min}, p_{\max})$, $\bar{c}=\bar{c}(p_{\min},
p_{\max}, c_1)$, $\varepsilon_0=\varepsilon_0(n,p_{\min}, p_{\max}, \bar{r}%
_1, c_1)$,\newline
$\varepsilon_1=\varepsilon_1(n,p_{\min}, p_{\max}, \bar{r}_1, c_0, c_1,
\theta)$.
\end{lemma}



\smallskip

Another essential tool leading to the proof of Theorem \ref{flatmain1fl} is given by the following Harnack type inequality, obtained in 
\cite{FL}, whose $p-$Laplace version (i.e., for $p(x)\equiv p$ constant) and  with $p\ge 2$, is contained in 
\cite{LR}.

\begin{lemma}
\label{harnack-with-e} Assume that $1<p_{\min}\le p(x)\le p_{\max}<\infty$
with $p(x)$ Lipschitz continuous in $\Omega$ and $\|\nabla
p\|_{L^{\infty}}\leq L$, for some $L>0$. Let $x_0\in\Omega$ and $0<R\le 1$
such that $\overline{B_{4R}(x_0)}\subset\Omega$. Let $v\in
W^{1,p(\cdot)}(\Omega)\cap L^{\infty}(\Omega)$ be a nonnegative solution to 
\begin{equation}  \label{eq-with-e}
\mbox{\rm div} (|\nabla v+e|^{p(x)-2}(\nabla v+e)) =f \quad\mbox{
in }\Omega,
\end{equation}
where $f\in L^{\infty}(\Omega)$ with $||f||_{L^{\infty}(\Omega)}\le 1$ and $%
e\in \mathbb{R}^n$ with $|e|=1$. Then, there exists $C$ such that 
\begin{equation}  \label{quasi_Harnack}
\sup_{{B_R}(x_0)}v\leq C\Big[\inf_{B_R(x_0)}v+ R\Big({||f||_{L^{%
\infty}(B_{4R}(x_0))}}^{\frac{1}{p_{\max}-1}}+C\Big)\Big].
\end{equation}
The constant $C$ depends only on $n$, $p_{\min}$, $p_{\max}$, $%
||v||_{L^{\infty}(B_{4R}(x_0))}$ and $L$.
\end{lemma}

Let us also mention that the fact that weak solutions to the inhomogeneous $p(x)$-Laplacian are locally of
class $C^{1,\alpha }$ plays a critical role in the results in \cite{FL}. We
point out that sharp conditions about the regularity of solutions of some
elliptic equations with non-standard growth can be found in \cite{AM} and 
\cite{Fan}.

\bigskip

Finally,  the overall strategy of the proof of the Theorem \ref{flatmain1fl} in 
\cite{FL} follows closely the one in \cite{D}
and, after a suitable rescaling, it relies on the following

\begin{lemma}[Improvement of flatness]
\label{improv1}Let $u$ satisfy \eqref{fb} in $B_1$ and 
\begin{equation}\label{ep-bound}
 \|f\|_{L^\infty(B_1)} \leq {\varepsilon}^2, \quad ||g-1||_{L^{\infty}(B_1)}\le {\varepsilon}^2, \quad  
||\nabla p||_{L^{\infty}(B_1)}\le {\varepsilon}^{1+\theta},
\quad ||p-p_0||_{L^{\infty}(B_1)}\le {\varepsilon},
\end{equation}
for $0<\varepsilon<1$, for some constant $0<\theta\le 1$. Suppose that 
\begin{equation}  \label{flat_1}
(x_n -\varepsilon)^+ \leq u(x) \leq (x_n + \varepsilon)^+ \quad \text{in $%
B_1,$} \quad 0\in F(u).
\end{equation}
If $0<r \leq r_0$ for $r_0$ universal, and $0<\varepsilon \leq \varepsilon_0$
for some $\varepsilon_0$ depending on $r$, then 
\begin{equation}  \label{improvedflat_2_new}
(x \cdot \nu -r \varepsilon / 2)^+ \leq u(x) \leq (x \cdot \nu +r
\varepsilon/2)^+ \quad \text{in $B_r,$}
\end{equation}
with $|\nu|=1$ and $|\nu - e_n| \leq \tilde C\varepsilon$ for a universal
constant $\tilde C.$
\end{lemma}

We schematize below the main steps of its proof.

\textit{Step 1: Compactness.} Fix $r \leq r_0$ with $r_0$ universal (chosen in Step 3). Assume
by contradiction that there exists a sequence $\varepsilon_k \rightarrow 0$
and a sequence $u_k$ of solutions to \eqref{fb} in $B_1$ with right hand
side $f_k$, exponent $p_k$ and free boundary condition $g_k$ satisfying %
\eqref{ep-bound} with $\varepsilon=\varepsilon_k$, such that $u_k$ satisfies %
\eqref{flat_1}, i.e., 
\begin{equation}  \label{flat_k}
(x_n -\varepsilon_k)^+ \leq u_k(x) \leq (x_n + \varepsilon_k)^+ \quad \text{%
for $x \in B_1$, $0 \in F(u_k),$}
\end{equation}
but $u_k$ does not satisfy the conclusion \eqref{improvedflat_2_new} of the
lemma.

Set 
\begin{equation*}
\tilde{u}_{k}(x)= \dfrac{u_k(x) - x_n}{\varepsilon_k}, \quad x \in
\Omega_1(u_k).
\end{equation*}
Then, \eqref{flat_k} gives 
\begin{equation}  \label{flat_tilde**}
-1 \leq \tilde{u}_{k}(x) \leq 1 \quad \text{for $x \in\Omega_1(u_k)$}.
\end{equation}

With a compactness argument and a sharp application of Ascoli-Arzel\`a
theorem it is possible to prove that there exists a convergent subsequence
to a function $\tilde{u}.$

\textit{Step 2. Limiting solution.} A delicate argument, involving our assumption \eqref{ep-bound}, shows that the function $\tilde u$ solves, in the viscosity sense, the
following linearized problem 
\begin{equation}  \label{Neumann_p-in-1/2}
\begin{cases}
{\mathcal{L}}_{p_0} \tilde u=0 & \text{in $B_{1/2} \cap \{x_n > 0\}$}, \\ 
\tilde u_n=0 & \text{on $B_{1/2} \cap \{x_n =0\}.$}%
\end{cases}%
\end{equation}

Here $1<p_{\min}\le p_0\le p_{\max}<\infty$, $\tilde u_n$ denotes the
derivative in the $e_n$ direction of $\tilde u$ and 
\begin{equation}  \label{Neumann_p-in-1/2-1/3}
{\mathcal{L}}_{p_0} u := \Delta u +(p_0-2) {\partial}_{nn} u.
\end{equation}
We point out that, if  
$\tilde u$ is a
continuous function on $B_{1/2}\cap\{x_n\ge 0\}$, we say that $\tilde u$ is
a viscosity solution to \eqref{Neumann_p-in-1/2}, if given a quadratic
polynomial $P(x)$ touching $\tilde u$ from below (resp. above) at $\bar{x}%
\in B_{1/2} \cap \{x_n \ge 0\}$,

\begin{itemize}
\item[(i)] if $\bar{x}\in B_{1/2} \cap \{x_n > 0\}$ then ${\mathcal{L}}%
_{p_0}P\le 0$ (resp. ${\mathcal{L}}_{p_0}P\ge 0$), i.e. ${\mathcal{L}}%
_{p_0}\tilde u= 0$ in the viscosity sense in $B_{1/2} \cap \{x_n > 0\}$;

\item[(ii)] if $\bar{x}\in B_{1/2} \cap \{x_n = 0\}$ then $P_n(\bar{x})\le 0$
(resp. $P_n(\bar{x})\ge 0$).
\end{itemize}


\textit{Step 3: Improvement of flatness.} From the previous step, $\tilde u$
solves \eqref{Neumann_p-in-1/2} and from \eqref{flat_tilde**}, 
\begin{equation}\label{bound-tilde-u}
-1 \leq \tilde{u}(x) \leq 1 \quad \text{in $B_{1/2}\cap\{x_n\ge 0\}$}.
\end{equation}

On the other hand 
it is well known that if $\tilde{u}$ is a viscosity solution to %
\eqref{Neumann_p-in-1/2} in $B_{1/2}\cap \{x_{n}\geq 0\}$, then $\tilde{u}%
\in C^{2}(B_{1/2}\cap \{x_{n}\geq 0\})$ and it is a classical solution to %
\eqref{Neumann_p-in-1/2}. 
Then, using the bound \eqref{bound-tilde-u}
 we find that, for the given $r$, 
\begin{equation*}
|\tilde{u}(x)-\tilde{u}(0)-\nabla \tilde{u}(0)\cdot x|\leq C_{0}r^{2}\quad 
\text{in }B_{r}\cap \{x_{n}\geq 0\}
\end{equation*}%
if $r_0\le 1/4$, with $C_0$ universal. Now the iterative argument is the same that has been applied in the
linear case, \cite{D}, \cite{DFS}.


\medskip 
Received xxxx 20xx; revised xxxx 20xx. \medskip

\end{document}